\documentclass[A4paper,12pt]{article}
\usepackage{latexsym}
\usepackage{mathrsfs}
\usepackage{amssymb}
\usepackage{amscd}
\usepackage[dvips]{graphicx}                  

\newtheorem{theorem}{Theorem}
\newtheorem{corollary}{Corollary}

\newenvironment{proof}[1][Proof]{\textbf{#1.} }{\ \rule{0.5em}{0.5em}}
\date{}
\long\def\symbolfootnote[#1]#2{\begingroup%
	\def\thefootnote{$\;$}\footnote[#1]{$^*$#2}\endgroup}
\begin{document}
	
	\title{Remarks on some special points in extremally disconnected compact spaces}
	\author{Joanna Jureczko}
\maketitle

\symbolfootnote[2]{Mathematics Subject Classification: 54D35, 54E45, 54C10.

	\hspace{0.2cm}
	}

\begin{abstract}
	 We prove under ZFC that in each extremally disconnected compact space there exists a non-limit point of any countable discrete subset. 
 	\end{abstract}

	\section{Introduction}
	This paper enlarges the paper \cite{RF} in which the following results under MA were obtained
	\begin{enumerate}
		\item In each extremally disconnected compact space there exists  a non-limit point of any countable discrete subset.
		\item In $G(I)$ (Gleason space over the interval $[0,1]$) there exists a non-limit point of any discrete subset of cardinality $<2^\omega$.
		\item Kunen's Theorem on non-limit points from \cite{KK1}. 
	\end{enumerate} 
In this paper we show that the first statement is true in ZFC, (without any additional assumptions). The second one can be obtained in ZFC in similar way, that is why we do not show it here, but the third one we omit, because of Kunen's result in \cite{KK}.

All the above results touch the problem of the  investigation of the non-homogeneity of $\omega^*$. We understand homogeneity in the following way: a space $X$ is \textit{homogenous} if for any distinct points $x, y \in X$ there is an automorphism $f \colon X \to Y$ such that $f(x) =y$.
	
	In the literature there are known  results concerning the existence of a particular type of points which make $\omega^*$ inhomogeneous.
	For example $P$-points are such points. (We recall that $p\in X$ is a \textit{$P(\lambda)$-point} iff $p \in int(\bigcap\mathcal{F}$) for any family $\mathcal{F}$ of cardinality less than $\lambda$ of neighourhoods of $p$. If $\lambda=\omega_1$, then we call it a $P$-point). 
	 In fact, a homeomorphism $f \colon X \to X$ maps $P$-points  into $P$-points.	W. Rudin, \cite{WR}, proved  (under CH) that there are $2^{2^\omega}$ $P$-points in $\omega^*$, but S. Shelah, \cite{SS},  constructed a model in ZFC in which $\omega^*$ has no $P$-points. Thus, there is impossible to prove in set theory whether $P$-points exist in $\omega^*$.
	 Fortunately, one can prove (without any additinal assumptions) that in $\omega^*$ there are weak $P$-points. (A point $p \in X$ is a \textit{weak $P$-point} iff is not a limit point of any countable subset of $X \setminus \{p\}$). Obviously, every $P$-point is a weak $P$-point, but not conversely, which was proved  by Kunen in \cite{KK}. For showing this, Kunen used a special kind of family which he called independent linked family and the notion of so called $OK$-points. 
	 
	In the main results of this paper,  there is used a modification of   
 Kunen's proof of \cite[Theorem 3.1]{KK} and some ideas from \cite{RF}.
	 
	 The paper is divided into three sections. In Section 2 there are gathered main definitions and facts used in the paper. Section 3 contains main result. Section 4 contains corollaries of the main result.
	 
	 All notations used in the paper are standard for the area and are assumed as well known. For definitions and facts not cited here the reader is referred to~\cite{RE, FZ, TJ}.   
	 
 	\section{Definitions and previous results}
 	\textbf{2.1.} Le $\kappa$ be a cardinal. Let $Y$ be a topological space and let $p \in Y$. A sequence $(U_n\colon n \in \omega)$ of neighbourhoods of $p$ is a \textit{$\kappa$-$OK$ sequence} iff there are neigbbourhoods $\{V_\alpha \colon \alpha<\kappa \}$ of $p$ such that  for all $n\geqslant 1$ and  all $\alpha_1 < \alpha_2 < ...< \alpha_n < \kappa$ we have
 	$$V_{\alpha_1}\cap V_{\alpha_2}\cap ... \cap V_{\alpha_n} \subset U_n.$$
 	
 	A point $p \in Y$ is an \textit{$\kappa$-$OK$-point} in $Y$ iff every sequence of $\omega$ many neighbourhoods of $p$ is $\kappa$-$OK$. 
 	\\
 	\\
 	\textbf{Fact 1 \cite{KK}.} \textit{There is a $p \in \omega^*$ which is $2^\omega$-$OK$.}
 	\\\\
 	\textbf{Fact 2 \cite{KK}.} \textit{Let $Y$ be a $T_1$-space and $p\in Y$. If $p$ is $2^\omega$-$OK$ in  $Y$, then $p$ is a weak $P$-point.}
 	\\
 	\\
 	\textbf{2.2.}
 	Consider  a family of infinite subsets of $\omega$
 	$$A = \{A_k(\nu, \eta) \colon \nu< 2^\omega, \eta < 2^\omega, k < \omega\},$$
 	For any $t \in 2^\omega \times 2^\omega$ consider $$A(t) = \bigcap_{\nu \in r(t)}\bigcap_{\eta \in t_\nu}A_{k_\nu}(\nu, \eta),$$
 	 where $r(t) = \{\nu < 2^\omega \colon \exists_\eta (\nu, \eta) \in t\}, t_\nu = \{\eta < 2^\omega \colon (\nu, \eta) \in t\}$ and $k_\nu$ is the number of elements of $t_\nu$. 
 	 
 	 Let $\mathcal{F}$ be a filter on $\omega$  and \textit{fin} be the ideal of all finite subsets of $\omega$. The family $A$ is \textit{independent linked with respect to $\mathcal{F}$} iff the following conditions are fulfilled
 	 \begin{enumerate}
 	 \item $A(t) \cap B$ is infinite for any $t \in  2^\omega \times 2^\omega$  and $B \in \mathcal{F}$,
 	 \item $A_k(\nu, \eta) \subseteq A_n(\nu, \eta)$ for any $k \leqslant n$ and  $\nu, \eta < 2^\omega$,
 	 \item $|\bigcap_{\eta \in t_\nu} A_m(\nu, \eta)| < \omega$ for   any $\nu < 2^\omega$ and $m <k_\nu$.
 	 \end{enumerate} 
 	 \textbf{Fact 3 \cite{FZ}.} \textit{There is a family $A$ which is independent linked with respect to fin}.

	\section{Main result}
	
	This section contains four results. The first two concern the considerations of the existing of $2^\omega$-$OK$-points for one function, while the last two concern the considerations of existing $2^\omega$-$OK$-points which are common for $2^\omega$ many functions. 
	Since proofs of these results are essentially different one decided to present all of them with all details.
	\\

 Let $X$ be a space, (since now, we assume that $X$ is a $T_2$-space), and let $$\{f_\alpha \colon f_\alpha \colon X \stackrel{onto}{\longrightarrow} \beta \omega, \alpha < 2^\omega\}$$ be a family of open mappings.

Let $\{S_\beta \colon \beta < 2^\omega\}$ be a family of non-empty and closed subsets of $X$. For each $\gamma < 2^\omega$ let 
$$C_\gamma = \bigcap\{S_\beta \colon \beta < \gamma\}$$
and
$$B_\gamma = \left\{\begin{array}{rcl}
C_\gamma \cap f^{-1}_\nu(\omega^*) & \textrm{where } \nu = \min \{\beta < 2^\omega \colon C_\gamma \cap f^{-1}_\beta(\omega^*) \not = \emptyset\}\\
C_\gamma & \textrm{otherwise}
\end{array} \right.$$

In Theorem 1 we will consider one function so we have to slightly modify the definition of $B_\gamma$, i.e. $B_\gamma = C_\gamma \cap f^{-1}(\omega^*)$ whenever $C_\gamma \cap f^{-1}(\omega^*) \not = \emptyset$.

\begin{theorem}
	Let $X$ be a compact space and let $ f \colon X \stackrel{onto}{\longrightarrow} \beta \omega$ be an open mapping. Then  there exists a family $\{U_{\beta} \colon \beta < 2^\omega\}$ of non-empty and closed subsets of $X$ such that for each $\gamma < 2^\omega$
	\begin{itemize}
		\item [(i)] $U_{\beta} \subseteq U_{\beta'}$ for any $\beta' \leqslant \beta < \gamma$;
		\item [(ii)] if $f^{-1}(\omega^*) \cap B_\gamma = \emptyset$, then
		$U_{\gamma} = f^{-1}(\{\xi\}) \cap B_\gamma$ for some $\xi \in \omega$;
		\item [(iii)] if $f^{-1}(\omega^*) \cap B_\gamma  \not = \emptyset$, then there exists an ultrafilter $\bigcup_{\beta \in 2^\omega}\mathcal{F}_\beta$ such that 
		\begin{itemize}
			\item [(a)] $\mathcal{F}_{\beta} \subseteq \mathcal{F}_{\beta'}$, whenever $\beta \leqslant \beta'< \gamma$
		\item [(b)] for a decreasing sequence $(V_\beta^k \colon k < \omega)$ if each $V^k_\beta \in \mathcal{F}_\beta$ then there exists $W^\eta_\beta \in \mathcal{F}_{\beta+1}$ for each $\eta \in 2^\omega$ such that for each $k \geqslant 1$ and each $\eta_1 < \eta_2 < ... < \eta_k < 2^\omega$ 
		$$|\bigcap_{i=1}^{k}W^{\eta_i}_{\beta} \setminus V^{k}_{\beta}|< \omega.$$
			\item [(c)] $f^{-1}(U) \cap B_\gamma \not =\emptyset$, where $U$ is a clopen set belonging to the standard base of $\omega^*$ such that $\bigcup_{\beta< \gamma} \mathcal{F}_\beta \in U$.
		\end{itemize}	
	\end{itemize}
\end{theorem}

		\begin{proof}	
	 Enumerate all subsets $\{Z_\alpha \colon \alpha < 2^\omega\}$  of $\omega$. Let $$\{(V_\beta^k \colon k < \omega) \colon \beta< 2^\omega\}$$ be a family  of  decreasing sequences, (i.e. $V^{k+1}_\beta \subset V^k_\beta$ for any $k < \omega$) and such that each decreasing sequence of subsets of $\omega$ is listed cofinally often. Let $$\{A_k(\nu, \eta) \colon k < \omega, \nu\in 2^\omega, \eta \in 2^\omega\}$$ be an independent linked family with respect to  fin. (We do not require the $(V^k_\beta\colon k<\omega)$ are to be distinct).

We construct $\mathcal{F}_\beta$, sets of indices $I_{\beta}$ and $U_{\beta}$ as follows, ($B_\beta$ are defined as is given before the theorem),

	\begin{itemize}
	\item[(1)] $\mathcal{F}_0$ is the cofinite filter, $I_0= 2^\omega$ and $U_{0} = X$,
	\item[(2)] $\mathcal{F}_\beta$ is a filter on $\omega$, $I_\beta \subset 2^\omega$,  $U_{\beta}$ is a non-empty closed subset of $X$ and $\{A_k(\nu, \eta) \colon \nu \in I_\beta, 1\leqslant k<\omega, \eta \in 2^\omega\}$ is an independent linked family with respect to $\mathcal{F}_\beta$
	\item[(3)] if $\gamma < \beta$, then $\mathcal{F}_\gamma \subseteq \mathcal{F}_\beta$, $I_\gamma\supseteq I_\beta$, $U_{\gamma} \supseteq U_{\beta}$,
	\item[(4)] if $\beta$ is limit, then $\mathcal{F}_\beta = \bigcup_{\gamma< \beta} \mathcal{F}_\gamma$, $I_\beta = \bigcap_{\gamma<\beta} I_\gamma$, $U_{\beta} = \bigcap_{\gamma<\beta} U_{\gamma}$,
	\item[(5)] $I_\beta \setminus I_{\beta+1}$ is finite for each $\beta \in 2^\omega$,
	\item[(6)] either $Z_\beta \in \mathcal{F}_\beta$ or $\omega\setminus Z_\beta \in \mathcal{F}_\beta$, 
\item[(7)] if $f(B_{\beta}) \cap \omega^* = \emptyset$ then $U_{\beta} = B_{\beta} \cap f^{-1}(\{\xi\})$ for some $\xi \in \omega$,
\item[(8)] if $f(B_{\beta}) \cap \omega^* \not = \emptyset$ then 
\begin{itemize}
	\item [(a)]  $f(B_{\beta}) \cap U \not =\emptyset$, for some clopen set $U$ belonging to the standard base of $\omega^*$ and $\bigcup_{\gamma < \beta} \mathcal{F}_\gamma \in U$,
	\item [(b)] $U_{\beta} = f^{-1}(U) \subseteq B_{\beta}$,
		\item [(c)] for a decreasing sequence $(V_\beta^k \colon k < \omega)$ if each $V^k_\beta \in \mathcal{F}_\beta$ then there exists $W^\eta_\beta \in \mathcal{F}_{\beta+1}$ for each $\eta \in 2^\omega$ such that for each $k \geqslant 1$ and each $\eta_1 < \eta_2 < ... < \eta_k < 2^\omega$ 
			$$|\bigcap_{i=1}^{k}W^{\eta_i}_{\beta} \setminus V^{k}_{\beta}|< \omega.$$
	\end{itemize}
\end{itemize} 	

The first and the limit step are done. Assume that we have constructed $\mathcal{F}_\beta$, $I_\beta$ and $U_{\beta}$. We will show the successor step.  
Assume that  $$f(B_{\beta}) \cap \omega^* \not = \emptyset.$$ Otherwise, by $(7)$ we put $U_{\beta+1} = B_{\beta+1} \cap f^{-1}(\{\xi\})$ for some $\xi \in \omega$. Then $(ii)$ holds.
\\
Consider two cases.

Case 1. $\beta \equiv 0 \ (mod\  2)$. Let $\mathcal{R}$ be the filter generated by $\mathcal{F}_\beta$ and $Z_\beta$, (abbr. $\mathcal{R} = [\mathcal{F}_\beta, \{Z_\beta\}]$).
If $\mathcal{R}$ is a proper filter and 
$$\{A_k(\nu, \eta) \colon \nu \in I_\beta, 1\leqslant k<\omega, \eta \in 2^\omega\}$$	
	is an independent linked family with respect to  $\mathcal{R}$ then we put
	$$\mathcal{F}_{\beta+1} = \mathcal{F}_\beta, \ \ I_{\beta+1} = I_\beta, \ \ U_{\beta+1} = U_{\beta}.$$
	If not, then there exists $E \in \mathcal{F}_\beta$ such that 
	$$Z_\beta\cap E \cap \bigcap_{\nu \in \tau}(\bigcap_{\eta \in \sigma_\nu} A_{k_\nu}(\nu, \eta)) = \emptyset$$
	for some $\tau \in [I_\beta]^{<\omega}, $ $ k_\nu \in \omega,$ and $ \sigma_\nu \in [2^\omega]^{k_\nu}$.	
	Then we put
$$\mathcal{F}_{\beta+1} = [\mathcal{F}_\beta, \bigcap_{\nu \in \tau}(\bigcap_{\eta \in \sigma_\nu} A_{k_\nu}(\nu, \eta))],$$ $$I_{\beta+1} = I_{\beta}\setminus~\tau$$
and
$$U_{\beta +1} = f^{-1}(U) \subseteq B_{\beta+1},$$
where $U$ is a clopen set belonging to the standard base of $\omega^*$ and $\bigcup_{\gamma < \beta} \mathcal{F}_\gamma \in U$.

Hence $\omega\setminus Z_\beta \in \mathcal{F}_{\beta+1}.$

Case 2. $\beta \equiv 1\ (mod\ 2)$.
If $V^k_\beta \not \in \mathcal{F}_\beta$ for some $k$, ($1\leqslant k < \omega$) then put 
	$$\mathcal{F}_{\beta+1} = \mathcal{F}_\beta, \ \ I_{\beta+1} = I_\beta, \ \ U_{\beta+1} = U_{\beta}.$$
	If $V^k_\beta \in \mathcal{F}_\beta$ for any $k$, ($1\leqslant k < \omega$) then fix $\nu \in I_\beta$, (by $(5)$ we have that $I_\beta$ is non-empty).
	Consider a set
		$$W_\beta^\eta = (\bigcap_{k} V_\beta^k)\cap(\bigcup_{1\leqslant k < \omega}(A_{k}(\nu, \eta)\cap V_\beta^k \setminus V_\beta^{k+1}))$$
		and put 
		$$\mathcal{F}_{\beta+1} = [\mathcal{F}_\beta, W^\eta_\beta],\ \ I_{\beta+1} = I_\beta \setminus \{\nu\}$$
		and $$U_{\beta+1} = f^{-1}(U) \subseteq B_{\beta+1},$$
		where $U$ is a clopen set belonging to the standard base of $\omega^*$ and $\bigcup_{\gamma < \beta} \mathcal{F}_\gamma \in U$,
		
		Now, we show that $(8)$ holds.
		Take $\eta_1< \eta_2<...< \eta_k< 2^\omega$ and consider 
		$$D = \bigcap_{i=1}^{k}W^{\eta_i}_{\beta} \setminus V^{k}_{\beta}.$$
		Obviously, for $k=1$ the set $D$ is empty. For $k>1$, we have 
		$$D \subseteq \bigcap_{i=1}^{k} A_{k-1}(\nu, \eta_i)$$
		is finite, (by the definition of an indepedent linked family).
		Moreover, 
		$$W^\eta_{\beta} \supseteq V^k_\beta \cap A_k(\nu, \eta)$$
		for each $k$, $(1\leqslant k < \omega)$ which verifies $(8)$.
\end{proof}			
\\			
	
\begin{theorem}
	Let $X$ be a compact space and let $f \colon X \stackrel{onto}{\longrightarrow} \beta \omega$ be an open mapping. Then there exists a family of sets $\{U_\beta \colon \beta < 2^\omega\}$  fulfilling $(i)-(iii)$ of Theorem 1 such that for each  $x \in \bigcap\{U_\beta \colon \beta < 2^\omega\}$  either $f(x)$ belongs to $\omega$ or $f(x)$ is $2^\omega$-$OK$-point in $\omega^*$.
\end{theorem}

\begin{proof}
	By Theorem 1, there exists a family $\{U_{\beta} \colon \beta \in 2^\omega\}$ of non-empty and closed subsets of $X$ and an ultrafilter $\bigcup_{\beta \in 2^\omega}\mathcal{F}_\beta$ of properties $(i)-(iii)$. 
	Then for each $x \in \bigcap\{U_{\beta} \colon \beta \in 2^\omega\}$ the image $f(x)$ belongs to $\omega$ whenever
	$f_{\alpha}(B_{\beta}) \cap \omega^* = \emptyset$ or $f(x)$    is a $2^\omega$-$OK$-point in $\omega^*$.
	(of the form $\bigcup_{\beta \in 2^\omega}\mathcal{F}_{\beta}$) whenever 
	$ f(B_{\beta}) \cap \omega^* \not = \emptyset$.
\end{proof}
\\

Let $<_*$ denote the canonical well-ordering on $2^\omega \times 2^\omega$. Let $j(\alpha, \beta)\in 2^\omega\times 2^\omega$ be an immediate successor of $(\alpha, \beta)$, i.e. $(\alpha, \beta)<_*j(\alpha, \beta)$ and there is no $(\gamma, \delta)$ such that $(\alpha, \beta)<_*(\gamma, \delta)<_*j(\alpha, \beta)$.

Let $X$ be a space and let $\{f_\alpha \colon f_\alpha \colon X \stackrel{onto}{\longrightarrow} \beta \omega, \alpha < 2^\omega\}$ be a family of open mappings. 
Let $\{S_{\alpha, \beta} \colon (\alpha, \beta) < 2^\omega\times 2^\omega\}$ be a family of non-empty and closed subsets of $X$. For each $(\gamma, \delta) \in 2^\omega \times 2^\omega$ let 
$$C_{\gamma, \delta} = \bigcap\{S_{\alpha,\beta} \colon (\alpha, \beta) <_* (\gamma, \delta)\}$$ and
$$B_{\gamma, \delta} = \left\{\begin{array}{rcl}
C_{\gamma, \delta} \cap f^{-1}_\mu(\omega^*) & \textrm{where } \mu = \min \{\beta < 2^\omega \colon C_{\gamma, \delta} \cap f^{-1}_\beta(\omega^*) \not = \emptyset\}\\
C_{\gamma, \delta} & \textrm{otherwise}
\end{array} \right.$$

\begin{theorem}
	Let $X$ be a compact space and let $\{f_\alpha \colon f_\alpha \colon X \stackrel{onto}{\longrightarrow} \beta \omega, \alpha < 2^\omega\}$ be a family of open mappings. Then  there exists a family $\{U_{\alpha,\beta} \colon (\alpha, \beta) < 2^\omega\times 2^\omega\}$ of non-empty and closed subsets of $X$ such that for each $(\gamma, \delta) < 2^\omega\times 2^\omega$
	\begin{itemize}
		\item [(i)] $U_{\alpha, \beta} \subseteq U_{\alpha', \beta'}$ for any $(\alpha',\beta') \leqslant_* (\alpha,\beta) <_* (\gamma, \delta)$;
		\item [(ii)] if $f^{-1}_\gamma(\omega^*) \cap B_{\gamma, \delta} = \emptyset$, then
		$U_{\gamma, \delta} = f^{-1}_\gamma(\{\xi\}) \cap B_{\gamma, \delta}$ for some $\xi \in \omega$;
		\item [(iii)] if $f^{-1}_\gamma(\omega^*) \cap B_{\gamma, \delta}  \not = \emptyset$, then there exists an ultrafilter $\bigcup_{(\alpha, \beta) \in 2^\omega\times 2^\omega}\mathcal{F}_{\alpha, \beta}$ such that 
		\begin{itemize}
			\item [(a)] $\mathcal{F}_{\alpha, \beta} \subseteq \mathcal{F}_{\alpha', \beta'}$, whenever $(\alpha,\beta) \leqslant_* (\alpha',\beta') <_* (\gamma, \delta)$;
		\item [(b)] for a decreasing sequence $(V_\beta^k \colon k < \omega)$ if each $V^k_\beta \in \mathcal{F}_\beta$ then there exists $W^\eta_{\alpha, \beta} \in \mathcal{F}_{j(\alpha, \beta)}$ for each $\eta \in 2^\omega$ such that for each $k \geqslant 1$ and each $\eta_1 < \eta_2 < ... < \eta_k < 2^\omega$ 
		$$|\bigcap_{i=1}^{k}W^{\eta_i}_{\alpha, \beta} \setminus V^{k}_{\alpha, \beta}|< \omega.$$
			\item [(c)] $f^{-1}(U) \cap B_{\gamma, \delta} \not =\emptyset$, where $U$ is a clopen set belonging to the standard base of $\omega^*$ such that $\bigcup_{(\alpha, \beta)<_* (\gamma, \delta)} \mathcal{F}_{\alpha, \beta} \in U$.
		\end{itemize}	
	\end{itemize}
\end{theorem}

\begin{proof}		
	Enumerate all subsets $\{Z_\alpha \colon \alpha < 2^\omega\}$  of $\omega$. Let $$\{(V_\beta^k \colon k < \omega) \colon \beta< 2^\omega\}$$ be a family  of  decreasing sequences, (i.e. $V^{k+1}_\beta \subset V^k_\beta$ for any $k < \omega$) and such that each decreasing sequence of subsets of $\omega$ is listed cofinally often. Let $$\{A_k(\nu, \eta) \colon k < \omega, \nu\in 2^\omega, \eta \in 2^\omega\}$$ be an independent linked family with respect to fin. (We do not require the $(V^k_\beta\colon k<\omega)$ are to be distinct).
	
	We construct $\mathcal{F}_{\alpha, \beta}$, sets of indices $I_{\alpha, \beta}$ and $U_{\alpha, \beta}$ as follows, ($B_{\alpha, \beta}$ is defined as given before the theorem),
	\begin{itemize}
		\item[(1)] $\mathcal{F}_{0,0}$ is the cofinite filter, $I_{0,0}= 2^\omega$ and $U_{0, 0} = X$,
	\item[(2)] $\mathcal{F}_{\alpha, \beta}$ is a filter on $\omega$, $I_{\alpha, \beta} \subset 2^\omega$,  $U_{\alpha, \beta}$ is a non-empty closed subset of $X$, and $\{A_k(\nu, \eta) \colon \nu \in I_{\alpha, \beta}, 1\leqslant k<\omega, \eta \in 2^\omega\}$ is an independent linked family with respect to. $\mathcal{F}_{\alpha, \beta}$,
	\item[(3)] if $(\gamma, \delta) <_* (\alpha,\beta)$, then $\mathcal{F}_{\gamma, \delta} \subseteq \mathcal{F}_{\alpha,\beta}$, $I_{\gamma, \delta}\supseteq I_{\alpha,\beta}$, $U_{\gamma, \delta} \supseteq U_{\alpha,\beta}$,
	\item[(4)] if $(\alpha,\beta)$ is limit, then $\mathcal{F}_{\alpha, \beta} = \bigcup_{(\gamma, \delta) <_* (\alpha,\beta)} \mathcal{F}_{\gamma, \delta}$, $I_{\alpha,\beta} = \bigcap_{(\gamma, \delta) <_* (\alpha,\beta)} I_{\gamma, \delta}$, $U_{\alpha, \beta} = \bigcap_{(\gamma, \delta) <_* (\alpha,\beta)} U_{\alpha, \beta}$,
	\item[(5)] $I_{\alpha, \beta} \setminus I_{j(\alpha, \beta)}$ is finite for each $(\alpha, \beta) \in 2^\omega\times 2^\omega$,
	\item[(6)] either $Z_{\alpha,\beta} \in \mathcal{F}_{\alpha, \beta}$ or $\omega\setminus Z_{\alpha, \beta} \in \mathcal{F}_{\alpha, \beta}$, 
	\item[(7)] if $f(B_{\alpha, \beta}) \cap \omega^* = \emptyset$ then $U_{\alpha, \beta} = B_{\alpha, \beta} \cap f_\alpha^{-1}(\{\xi\})$ for some $\xi \in \omega$,
	\item[(8)] if $f_\alpha(B_{\alpha, \beta}) \cap \omega^* \not = \emptyset$ then 
	\begin{itemize}
		\item [(a)] $f_\alpha(B_{\alpha, \beta}) \cap U \not =\emptyset$, for some clopen set $U$ belonging to the standard base of $\omega^*$ and $\bigcup_{(\gamma, \delta) <_* (\alpha, \beta)} \mathcal{F}_{\gamma, \delta} \in U$,
		\item [(b)] $U_{\alpha, \beta} = f_\alpha^{-1}(U) \subseteq B_{\alpha, \beta}$
		\item [(c)] for a decreasing sequence $(V_{\alpha, \beta}^k \colon k < \omega)$ if each $V^k_{\alpha, \beta} \in \mathcal{F}_{\alpha, \beta}$ then there exists $W^\eta_{\alpha, \beta} \in \mathcal{F}_{j(\alpha, \beta)}$ for each $\eta \in 2^\omega$ such that for each $k \geqslant 1$ and each $\eta_1 < \eta_2 < ... < \eta_k < 2^\omega$ 
		$$|\bigcap_{i=1}^{k}W^{\eta_i}_{\alpha, \beta} \setminus V^{k}_{\alpha, \beta}|< \omega.$$
	\end{itemize}
	\end{itemize} 	
	
	The first and the limit step are done. Assume that we have constructed $\mathcal{F}_{\alpha, \beta}$, $I_{\alpha, \beta}$ and $U_{\alpha, \beta}$. We will show the successor step.  
	
	Assume that  $f_\alpha(B_{\alpha, \beta}) \cap \omega^* \not = \emptyset$. Otherwise, by $(7)$ we put $$U_{j(\alpha, \beta)} = B_{j(\alpha, \beta)} \cap f^{-1}_\alpha(\{\xi\})$$ for some $\xi \in \omega$. Then $(ii)$ holds.
	
	Consider two cases.
	
	Case 1. $\alpha +\beta \equiv 0 \ (mod\  2)$. Let $\mathcal{R}$ be the filter generated by $\mathcal{F}_{\alpha, \beta}$ and $Z_{\alpha, \beta}$, (abbr. $\mathcal{R} = [\mathcal{F}_{\alpha, \beta}, \{Z_{\alpha, \beta}\}]$).
	If $\mathcal{R}$ is a proper filter and 
	$$\{A_k(\nu, \eta) \colon \nu \in I_{\alpha, \beta}, 1\leqslant k<\omega, \eta \in 2^\omega\}$$	
	is an independent linked family with respect to  $\mathcal{R}$ then we put
	$$\mathcal{F}_{j(\alpha, \beta)} = \mathcal{F}_{\alpha, \beta}, \ \ I_{j(\alpha, \beta)} = I_{\alpha, \beta}, \ \ U_{j(\alpha, \beta)} = U_{\alpha, \beta}.$$
	If not, then there exists $E \in \mathcal{F}_{\alpha, \beta}$ such that 
	$$Z_{\alpha, \beta}\cap E \cap \bigcap_{\nu \in \tau}(\bigcap_{\eta \in \sigma_\nu} A_{k_\nu}(\nu, \eta)) = \emptyset$$
	for some $\tau \in [I_\beta]^{<\omega}, $ $ k_\nu \in \omega,$ and $ \sigma_\nu \in [2^\omega]^{k_\nu}$.	
	Then we put
	$$\mathcal{F}_{j(\alpha, \beta)} = [\mathcal{F}_{\alpha, \beta}, \bigcap_{\nu \in \tau}(\bigcap_{\eta \in \sigma_\nu} A_{k_\nu}(\nu, \eta))],$$ $$I_{j(\alpha, \beta)} = I_{\alpha, \beta}\setminus~\tau$$
	and
	$$U_{j(\alpha, \beta)} = f_\alpha^{-1}(U) \subseteq B_{j(\alpha, \beta)},$$
	where $U$ is a clopen set belonging to the standard base of $\omega^*$ and $\bigcup_{(\gamma, \delta) <_* (\alpha, \beta)} \mathcal{F}_{\gamma, \delta} \in U$,
	
	Hence $\omega\setminus Z_{\alpha, \beta} \in \mathcal{F}_{j(\alpha, \beta)}.$
	
	Case 2. $\alpha+\beta \equiv 1\ (mod\ 2)$.
	If $V^k_{\alpha, \beta} \not \in \mathcal{F}_{\alpha, \beta}$ for some $k$, ($1\leqslant k < \omega$) then put 
	$$\mathcal{F}_{j(\alpha, \beta)} = \mathcal{F}_{\alpha, \beta}, \ \ I_{j(\alpha, \beta)} = I_\beta, \ \ U_{j(\alpha, \beta)} = U_{\alpha, \beta}.$$
	If $V^k_{\alpha, \beta} \in \mathcal{F}_{\alpha, \beta}$ for any $k$, ($1\leqslant k < \omega$) then fix $\nu \in I_{\alpha, \beta}$, (by $(5)$ we have that $I_{\alpha, \beta}$ is non-empty).
	Consider a set
	$$W_{\alpha, \beta}^\eta = (\bigcap_{k} V_{\alpha, \beta}^k)\cap(\bigcup_{1\leqslant k < \omega}(A_{k}(\nu, \eta)\cap V_{\alpha, \beta}^k \setminus V_{\alpha, \beta}^{k+1}))$$
	and put 
	$$\mathcal{F}_{j(\alpha, \beta)} = [\mathcal{F}_{\alpha, \beta}, W^\eta_{\alpha, \beta}]\ \ I_{j(\alpha, \beta)} = I_{\alpha, \beta} \setminus \{\nu\}$$
	and $$U_{j(\alpha, \beta)} = f_\alpha^{-1}(U) \subseteq B_{j(\alpha, \beta)},$$
	where $U$ is a clopen set belonging to the standard base of $\omega^*$ and $\bigcup_{(\gamma, \delta) <_* (\alpha, \beta)} \mathcal{F}_{\gamma, \delta} \in U$,
	
	Now, we show that $(8)$ holds.
	Take $\eta_1< \eta_2<...< \eta_k< 2^\omega$ and consider 
	$$D = \bigcap_{i=1}^{k}W^{\eta_i}_{\alpha, \beta} \setminus V^{k}_{\beta}.$$
	Obviously, for $k=1$ the set $D$ is empty. For $k>1$, we have 
	$$D \subseteq \bigcap_{i=1}^{k} A_{k-1}(\nu, \eta_i)$$
	is finite, (by the definition of an indepedent linked family).
	Moreover, 
	$$W^\eta_{\alpha, \beta} \supseteq V^k_\beta \cap A_k(\nu, \eta)$$
	for each $k$, $(1\leqslant k < \omega)$ which verifies $(8)$.
\end{proof}
\\

\begin{theorem}
	Let $X$ be a compact space and let $$\{f_\alpha \colon \alpha< 2^\omega, f_\alpha \colon X \stackrel{onto}{\longrightarrow} \beta \omega\}$$ be a family of open mappings. Then there exists a family $\{U_{\alpha, \beta} \colon (\alpha, \beta) \in 2^\omega\times 2^\omega\}$ such that for each $x \in \bigcap \{U_{\alpha, \beta} \colon (\alpha, \beta) \in 2^\omega\times 2^\omega\}$ and for each $\alpha \in 2^\omega$  either $f_\alpha(x)$ belongs to $\omega$ or $f_\alpha(x)$ is  a$2^\omega$-$OK$-point in $\omega^*$.
\end{theorem}

\begin{proof}
	By Theorem 3, there exists a family $\{U_{\alpha, \beta} \colon (\alpha, \beta) \in 2^\omega\times 2^\omega\}$ of non-empty and closed subsets of $X$ and an ultrafilter $\bigcup_{(\alpha, \beta) \in 2^\omega\times 2^\omega}\mathcal{F}_{\alpha, \beta}$ of properties $(i)-(iii)$. 
	
	Then for each $x \in \bigcap\{U_{\alpha, \beta} \colon (\alpha, \beta) \in 2^\omega\times 2^\omega\}$
	and for each $\alpha \in 2^\omega$, $f_\alpha(x)$ belongs to $\omega$ whenever
	$f_{\alpha}(B_{\alpha, \beta}) \cap \omega^* = \emptyset$ or $f_\alpha(x)$    is a $2^\omega$-$OK$-point in $\omega^*$
	(of the form $\bigcup_{(\alpha, \beta) \in 2^\omega \times 2^\omega}\mathcal{F}_{\alpha, \beta}$) whenever 
	$ f_{\alpha}(B_{\alpha, \beta}) \cap \omega^* \not = \emptyset$.
\end{proof}
\\

\textbf{Remark.} Using the notion of $2^\omega$-soso points introduced in \cite{JK}, we can modify the proof of Theorem 1 and Theorem 3 (by changing (8)). 

Recall that a point $p \in Y$ is \textit{$2^\omega$-soso} in $Y$ iff for every countable family $\mathcal{W} $	of neighbourhoods of $p$, there is an $2^\omega$-$OK$ sequence $(U_n \colon n \in \omega)$ of neighbourhoods of $p$ such that $\mathcal{W} \subset (U_n \in n \in \omega)$.

Obviously $2^\omega$-$OK$ implies $2^\omega$-soso. Moreover, if $Y$ is a $T_1$ space and $p\in Y$  is $2^\omega$-soso, then $p$ is a weak $P$-point. However, the notion of $2^\omega$-soso point is more complicated then $2^\omega$-$OK$ point, hence the modification the proof of Theorem 1 and Theorem 3 are necessary.

		\section{Corollaries}
		
			Applying Theorem 4 and Fact 2 we have immediately the following result.
		
		\begin{corollary}
				Let $X$ be a compact space and let $$\{f_\alpha \colon \alpha< 2^\omega, f_\alpha \colon X \stackrel{onto}{\longrightarrow} \beta \omega\}$$ be a family of open mappings. Then there exists a point  $x \in X$ such that for each $\alpha < 2^\omega$ either $f_\alpha(x) = \xi$ for some $\xi \in \omega$ or $f_\alpha(x)$ is a weak $P$-point of $\omega^*$.
		\end{corollary}
		
		Applying Theorem 4 and modification to the proof of \cite[Theorem 3]{RF} we  obtain the following result.
		
		\begin{corollary}
			Let $X$ be a compact extremally disconnected space of weight $2^\omega$. There exists a point $x \in X$ such that $x$ is not a limit point of any countable discrete subset of $X$.
		\end{corollary}
	
	Let $U(m) = st(\mathcal{P}(m)/\mathcal{P}^{<m}(m))$, where $m$ is a cardinal and $\mathcal{P}^{<m}(m)$ is the ideal of all subsets of cardinality less than $m$. The space $U(m)$ is an $F$-space. If $m$ is regular, then $U(m)$ is embedded in $\beta m = st(\mathcal{P}(m))$ as a $\mathcal{P}(m)$-set.
	
	Applying Theorem 4 and modification to the proof of \cite[Theorem 4]{RF} we  obtain the following result.
	
	\begin{corollary}
		In the space $U(m)$, where $\omega < cf(m)\leqslant m < 2^\omega$, there exists a point $x$ which is not a limit point of any strongly discrete subset of $U(m)$ of cardinality not greater than $m$.   
	\end{corollary}

	\noindent
	{\sc Joanna Jureczko}
	\\
	Wroc\l{}aw University of Science and Technology, Wroc\l{}aw, Poland.
	\\
	{\sl e-mail: joanna.jureczko@pwr.edu.pl}
	
\end{document}